\input amstex
\documentstyle{amsppt}
\input xypic
\xyoption{all}

\def \ni{\noindent}

\def \ot{\otimes}
\def \ov{\overline}
\def \wt{\widetilde}
\def \wh{\widehat}

\def \de{\delta}
\def \ep{\epsilon}

\def \De{\Delta}

\def\Ad{\operatorname{Ad}}

\def \coev{\operatorname{coev}}

\def \ev{\operatorname{ev}}

\def \fin{\operatorname{fin}}
\def \ide{\operatorname{id}}

\def \Ima{\operatorname{Im}}

\def \Z{\operatorname{Z}}

\def \op{\operatorname{op}}
\def \cop{\operatorname{cop}}

\def \xcirc{\,}

\topmatter

\title
Semiquasitriangular Hopf algebras
\endtitle

\author
     Jorge A. Guccione and Juan J. Guccione
\endauthor

\address
     Jorge Alberto Guccione, Departamento de Matem\'atica, Facultad de
     Ciencias Exactas y Naturales, Pabell\'on 1 - Ciudad Universitaria,
     (1428) Buenos Aires, Argentina.
\endaddress

\email
     vander\@dm.uba.ar
\endemail

\address
     Juan Jos\'e Guccione, Departamento de Matem\'atica, Facultad de
     Ciencias Exactas y Naturales, Pabell\'on 1 - Ciudad Universitaria,
     (1428) Buenos Aires, Argentina.
\endaddress

\email
     jjgucci\@dm.uba.ar
\endemail

\abstract We say that a Hopf algebra $H$ is {\it
semicocommutative} if the right adjoint coaction factorizes
through $H\ot \Z(H)$, where $\Z(H)$ denotes the centre of $H$. For
instance the commutative and the cocommutative Hopf algebras are
semicocommutative. The quasitriangular Hopf algebras generalize
the cocommutative Hopf algebras. In this paper we introduce and
begin the study of a similar generalization for the
semicocommutative ones. These algebras, which we call {\it
semiquasitriangular} Hopf algebras have many of the basic
properties of the quasitriangular ones. In particular, they have
associated braided categories of representations in a natural way.
\endabstract

\subjclass\nofrills{{\rm 2000} {\it Mathematics Subject
Classification}.\usualspace} Primary 16W30; Secondary 18D10
\endsubjclass

\thanks
Supported by UBACYT X193 and CONICET
\endthanks

\keywords Hopf algebra\endkeywords

\endtopmatter

\document

\head Introduction \endhead
Let $k$ be a field. All the algebras and vector spaces considered
in this paper are over $k$, all the maps are $k$-linear maps and
the unadorned tensor product will denote the tensor product over
$k$. As usual, given a Hopf algebra $H$, we write $\mu$, $\eta$,
$\Delta$, $\epsilon$ and $S$, adorned with a subscript if
necessary, to denote the multiplication, the unit, the
comultiplication, the counit and the antipode of $H$,
respectively. For the comultiplication we use the Sweedler
notation $\De(h) = h_1\ot h_2$, without summation symbol.
Moreover, given a right $H$-comodule $M$ with coaction $\nu_M$, we
write $\nu_M(m) = m_0\ot m_1$ (also without any summation symbol).
If $M$ is a left $H$-module, then the notation that we will use is
$\nu_M(m)= m_{-1}\ot m_0$.

\smallskip

Recall that a Hopf algebra $H$ is called cocommutative if $\Delta
= \tau \xcirc \Delta$, where $\Delta$ denotes the comultiplication
of $H$ and $\tau\:H\ot H\to H\ot H$ is the flip $\tau(h\ot l)=l\ot
h$. Let $S$ be the antipode of $H$. It is easy to check that $H$
is cocommutative if and only if $\Ad(h):=h_2\ot S(h_1)h_3 = h\ot
1$ for all $h\in H$. That is, if the right adjoint coaction is
trivial. In fact, in this case, $S^2(h) = S^2(h_2)S(h_1)h_3 = h$
for all $h\in H$, and then $h_2\ot h_1 = h_2\ot S^2(h_1) = h_3\ot
S^2(h_2)S(h_1)h_4 = h_1\ot h_2$. The converse assertion is
trivial.

Let $H$ be a Hopf algebra and let $\Z(H)$ the centre of $H$. We
say that a Hopf algebra $H$ is {\it semicocommutative} if $\Ad(h)
\in H\ot \Z(H)$ for all $h\in H$. That is, if the right adjoint
coaction factorizes through $H\ot \Z(H)$.

\smallskip

For instance, the commutative and cocommutative Hopf algebras are
semicocommutative Hopf algebras. Moreover, the class of these
algebras is closed under the operations of taking tensor products,
subHopfalgebras and quotients.

\smallskip

Let $H$ be a semicocommutative Hopf algebra. A {\it normal
$H$-module} is a vector space $M$, endowed with a left action
$\rho\:H\ot M\to M$ and a right coaction $\nu\:M\to M\ot H$, such
that for all $h\in H$ and $m\in M$,

\smallskip

\roster

\item $\nu(m)\in M\ot \Z(H)$,

\smallskip

\item $\nu(h\cdot m)= h_2\cdot m_0\ot S(h_1)h_3m_1$.

\endroster

\smallskip

\noindent A {\it morphism $f\:M\to N$ of normal $H$-modules} is a
map which is a morphism of $H$-modules and $H$-comodules.

\smallskip

For instance, $H$, endowed with the left regular action and the
right adjoint coaction, is a normal $H$-module, and, if $H$ is a
cocommutative Hopf algebra, then each left $H$-module, endowed
with the trivial coaction, is a normal $H$-module.

Note that the definition of normal $H$-module is similar to the
definition of Yetter-Drinfeld module. In fact, when $H$ is
commutative both notions coincide.

\smallskip

Next, we mention without proof some results about
semicocommutative Hopf algebras.

\proclaim{Theorem} For each semicocommutative Hopf algebra $H$,
the category ${}_H\Cal M^{H}(n)$, of normal $H$-modules, is a
braided category. The unit object is $k$, endowed with the trivial
action and the trivial coaction, and the tensor product is the
usual tensor product over $k$, endowed with the diagonal action
and the diagonal coaction. The associative and unit constraints
are the usual ones and the braid $c$ is given by $c_{MN}(m\ot n)=
n_0\ot n_1\cdot m$.
\endproclaim

\example{Example} Let $G$ be a group and let $\Z(G)$ be the center
of $G$. A normal $k[G]$-module is a direct sum $M=\oplus_{g\in
\Z(G)} M_g$, of left $G$-modules $M_g$. A map $f\:M\to N$ is a
morphism of normal $k[G]$-modules if it is $G$-linear and
$f(M_g)\subseteq N_g$ for all $g\in \Z(G)$. The braid
$c_{MN}\:M\ot N\to N\ot M$ is given by $c_{MN}(m\ot n)= n\ot
g\cdot m$, for $m\in M$ and $n\in N_g$.
\endexample

\example{Example} Let  $G$ be a finite group. Using that a right
$k[G]^*$-comodule is the same that a left $k[G]$-module and that a
$k[G]^*$-module is the same that a $G$-graduate $k$-module, it is
easy to check that a normal $k[G]^*$-module is a left $G$-module
$M$, endowed with a decomposition $M=\oplus_{g\in G} M_g$, such
that
$$
m\in M_g \Rightarrow y\cdot m\in M_{ygy^{-1}}, \text{ for all
$g,y\in G$.}
$$
Moreover, a map $f\:M\to N$ is a morphism of normal $k[G]$-modules
if it is $G$-linear and $f(M_g)\subseteq N_g$ for all $g\in G$.
Finally, the braid $c_{MN}\:M\ot N\to N\ot M$ is given by
$c_{MN}(m\ot n)= g\cdot n \ot m$, for $m\in M_g$ and $n\in N$.

\endexample

\proclaim{Theorem} For each semicocommutative Hopf algebra $H$,
the full braided subcategory of ${}_H\Cal M^{H}(n)$, consisting of
all finite dimensional normal $H$-modules, is rigid.
\endproclaim

\proclaim{Proposition} If $H$ is an semicocommutative Hopf
algebra, then
$$
S^2(h)=h_2S(h_1)h_3,\qquad\text{for all $h\in H$.}
$$
\endproclaim

One can think the quasitriangular Hopf algebras \cite{D} as a
generalization of the cocommutative ones. The aim of this paper is
to introduce and begin the study of a similar generalization for
the semicocommutative ones, that we call semiquasitriangular Hopf
algebras. Such an algebra is a pair $(H,R)$, consisting of a Hopf
algebra $H$ and an invertible element $R$ of $H\ot H$, satisfying
suitable conditions. From the definition it follows that $(H,1\ot
1)$ is semiquasitriangular if and only if $H$ is
semicocommutative. Moreover, the quasitriangular Hopf algebras are
semiquasitriangular.

Our main results are Theorem~2.5, Corollary~2.10, Theorem~2.11 and
Proposition~3.2. In particular we get generalizations of the
semicocommutative Hopf algebras results mentioned above.

\smallskip

We want to note that one can arrive to the notion of
semiquasitriangular Hopf algebra in a different way to the
considered in this paper. In \cite{G-G1} (see also \cite{G-G2}) we
define a notion of Hopf crossed products, that generalize the
classical one introduced in \cite{B-C-M} and \cite{D-T}. In
\cite{M1} was proved that if $H$ admits a quasitriangular
structure, then the Drinfeld double $D(H)$ of $H$ is isomorphic to
a classical Hopf crossed product $A\# H$. In \cite{D-G-G} was
proved that for $D(H)$ be isomorphic to a Hopf crossed product
$A\# H$ in the sense of \cite{G-G1}, it suffices that $H$ admits a
semiquasitriangular structure. This gives a version for this
setting of the Majid result.

\head 1. Semiquasitriangular Hopf algebras \endhead
In this section we introduce the notion semiquasitriangular Hopf
algebras and study its basic properties. Moreover, we show that
this concept includes the ones of quasitriangular and
semicocommutative Hopf algebras. All the results of this section
are immediate in the last case.

\smallskip

Before beginning we establish some notations. Let $H$ be a Hopf
algebra, $R = \sum_i R_i^{(1)}\ot R_i^{(2)}$ an invertible element
of $H\ot H$ and $\tau\:H\ot H\to H\ot H$ the flip $\tau(h\ot l) =
l\ot h$.

\roster

\smallskip

\item We will write $R = R^{(1)}\ot R^{(2)}$, understanding
the summation symbol and the index $i$. Similarly $R^{-(1)}\ot
R^{-(2)}$ denotes $R^{-1}$. When it is necessary we let
$R^{'(1)}\ot R^{'(2)}$, $\ov{R}^{(1)}\ot \ov{R}^{2)}$, etcetera
denote copies of $R$.

\smallskip

\item Given $n\ge 1$ let $H^{\ot n}$ be the tensor products of
$n$ copies of $H$. For any $2$-tuple $(k_1,k_2)$ of distinct
elements of $\{1,\dots,n\}$, we let $R_{k_1,k_2}$ denote the
element of $H^{\ot n}$, given by
$$
R_{k_1,k_2} = \sum_i y_i^{(1)}\ot\cdots\ot y_i^{(n)},
$$
where $y_i^{(k_j)} = R_i^{(j)}$ for any $j = 1,2$ and $y_i^{(k)} =
1$, otherwise.

\smallskip

\item Given $1\le i < n$ we let $\tau_{in}\: H^{\ot n}\to H^{\ot
n}$ denote the map $\tau_{in}:= H^{\ot i-1}\ot \tau\ot H^{\ot
n-i-1}$. Here we extend the notation of $H^{\ot n}$, taking
$H^{\ot 0} = k$.

\endroster

\definition{Definition 1.1} A semiquasitriangular Hopf algebra is a pair
$(H,R)$, where $H$ is a Hopf algebra with bijective antipode and
$R\in H\ot H$ is an invertible element satisfying

\smallskip

\roster

\item $R^{(1)}_1\ot R^{(1)}_2\ot R^{(2)} = R^{(1)}\ot
R'{}^{(1)}\ot
R^{(2)}R'{}^{(2)}$,

\smallskip

\item $R^{(1)}\ot R^{(2)}_1\ot R^{(2)}_2  = R^{(1)}R'{}^{(1)}\ot
R'{}^{(2)}\ot R^{(2)}$,

\smallskip

\item $R^{(1)}\ot R^{(2)}_2 R'{}^{(1)}\ot R^{(2)}_1 R'{}^{(2)} =
R^{(1)}\ot R'{}^{(1)}R^{(2)}_1\ot R'{}^{(2)}R^{(2)}_2$,

\smallskip

\item $R^{(1)}_2 R'{}^{(1)}\ot R^{(1)}_1 R'{}^{(2)}\ot R^{(2)} =
R'{}^{(1)}R^{(1)}_1 \ot R'{}^{(2)} R^{(1)}_2\ot R^{(2)}$,

\smallskip

\item $\nu(h):= R^{(2)}h_2 R'{}^{(2)}\ot S(h_1) S(R^{(1)})
h_3 R'{}^{(1)} \in H\ot \Z(H)$ for all $h\in H$,

\smallskip

\item $\nu(h)= R^{(1)}h_2 R'{}^{(1)}\ot S(R'{}^{(2)})S(h_1)
R^{(2)}h_3$ for all $h\in H$
\endroster

\smallskip

\ni If $(H,R)$ is a semiquasitriangular Hopf algebra, we say that
$R$ is a semiquasitriangular structure for $H$.

\enddefinition

The following result generalizes items~(1)--(4) of the above
definition

\proclaim{Proposition 1.2} Let $(H,R)$ be a semiquasitriangular
Hopf algebra and let $n\ge 1$. The following assertions are valid:

\smallskip

\roster

\item $(\De^n\ot H)(R) = R_{1,n+2}R_{2,n+2} \cdots R_{n+1,n+2}$,

\smallskip

\item $(H\ot \De^n)(R) = R_{1,n+2}R_{1,n+1} \cdots R_{12}$,

\smallskip

\item $(H\ot \tau_{in}\xcirc \De^n)(R)R_{i+1,i+2} = R_{i+1,i+2}
(H\ot \De^n)(R)$ for all $1\le i\le n$,

\smallskip

\item $(\tau_{in}\xcirc \De^n\ot H)(R)R_{i,i+1} = R_{i,i+1}
(\De^n\ot H)(R)$ for all $1\le i\le n$,

\smallskip

\endroster

\endproclaim

\demo{Proof} Items~(1) and (2) follows easily by induction on $n$.
We prove item~(3) and leave the last one to the reader. Assume by
induction the formula is true for $n$. Since,
$$
(H\ot \De^{n-1})(R) = R_{1,n+1}R_{1n} \cdots R_{12},\tag 2
$$
we have
$$
\align (H & \ot \tau_{in}\xcirc \De^n)(R)R_{i+1,i+2} = (H^i\ot
\De^{\cop}
\ot H^{n-i})(R_{1,n+1}R_{1n} \cdots R_{12})R_{i+1,i+2}\\
& = R_{1,n+2}R_{1,n+1} \cdots R_{1,i+3}(H^i\ot \De^{\cop})
(R_{1,i+1})R_{i+1,i+2} R_{1i}R_{1,i-1} \cdots R_{12}.
\endalign
$$
Hence, by item~(3) of Definition~1.1,
$$
\align
(H &\ot \tau_{in}\xcirc \De^n)(R)R_{i+1,i+2} \\
& = R_{1,n+2}R_{1,n+1} \cdots R_{1,i+3}R_{i+1,i+2}(H^i\ot
\De)(R_{1,i+1})
R_{1i}R_{1,i-1} \cdots R_{12}\\
& = R_{i+1,i+2}(H^i\ot \De\ot H^{n-i})(R_{1,n+1}R_{1n} \cdots
R_{12})\\
& = R_{i+1,i+2}(H\ot \De^n)(R),
\endalign
$$
where the last equality follows by $(2)$.\qed
\enddemo

\proclaim{Proposition 1.3} If $(H,R)$ is a semiquasitriangular
Hopf algebra, then
$$
(\epsilon\ot H)(R)=(H\ot\epsilon)(R)=1, \quad (S\ot H)(R)
=R^{-1}\quad\text{and}\quad (H\ot S)(R^{-1})=R.
$$
Hence, $(S\ot S)(R)=R$.
\endproclaim

\demo{Proof} The proof given in \cite{M2, Lemma 2.1.2} for
quasitriangular Hopf algebras only use items~(1) and (2) of
Definition~1.1.\qed
\enddemo

\proclaim{Proposition 1.4}If $(H,R)$ is an semiquasitriangular
Hopf algebra, then
$$
R_{12}R_{13}R_{23}= R_{23}R_{13}R_{12}.
$$
\endproclaim

\demo{Proof}The proof given in \cite{M2, Lemma 2.1.4} for
quasitriangular Hopf algebras only use items~(2) and (3) of
Definition~1.1.\qed
\enddemo

\example{Example 1.5} A Hopf algebra $H$ is semicocommutative if
and only if $(H,1_{H\ot H})$ is a semiquasitriangular Hopf
algebra.
\endexample

\example{Example 1.6} Each quasitriangular Hopf algebra is
semiquasitriangular. In fact, if $(H,R)$ is a quasitriangular Hopf
algebra, then it is well known that $(H,R)$ satisfies conditions
(1)--(4) of Definition~1.1. Moreover, for all $h\in H$,
$$
\allowdisplaybreaks \align \nu(h)&=R^{(2)}h_2R'{}^{(2)}\ot S(h_1)
S(R^{(1)})
h_3R'{}^{(1)}\\
&= R^{(2)}h_2R'{}^{(2)}\ot  S(R^{(1)}h_1)
h_3R'{}^{(1)}\\
&= h_1R^{(2)}R'{}^{(2)}\ot  S(h_2R^{(1)})
h_3R'{}^{(1)}\\
&= h R^{(2)}R'{}^{(2)}\ot  S(R^{(1)})R'{}^{(1)}\\
&=h\ot 1,
\endalign
$$
which clearly belongs to $H\ot \Z(H)$, and
$$
\allowdisplaybreaks \align R^{(1)}h_2R'{}^{(1)}\ot
S(R'{}^{(2)})S(h_1) R^{(2)} h_3& =
h_3R^{(1)}R'{}^{(1)}\ot S(R'{}^{(2)})S(h_1)h_2R^{(2)}\\
&= hR^{(1)}R'{}^{(1)}\ot S(R'{}^{(2)})R^{(2)}\\
&= hR^{(1)}R'{}^{(1)}\ot S(S^{-1}(R^{(2)})R'{}^{(2)})\\
&=h\ot 1.
\endalign
$$
\endexample

\example{Example 1.7} If $(H,R_H)$ and $(L,R_L)$ are
semiquasitriangular Hopf algebras, then $(H\ot L, R_H\ov{\ot}
R_L)$, where $R_H\ov{\ot} R_L:= (H\ot \tau\ot L)(R_H\ot R_L)$,
also is.
\endexample

\example{Example 1.8} The class of semiquasitriangular Hopf
algebra is closed under the operations of taking quotients.
Moreover, if $(H,R)$ is a quasitriangular Hopf algebra,
$L\subseteq H$ is a subHopfalgebra of $H$ and $R\in L\ot L$, then
$(L,R)$ is a semiquasitriangular Hopf algebra.
\endexample

\proclaim{Proposition~1.9} If $(H,R)$ is a semiquasitriangular
Hopf algebra, then
$$
\align \nu(h)&= R^{(2)}h_2R'{}^{(2)}\ot S^{2}(R'{}^{(1)})S(h_1)
S(R^{(1)})h_3\\
&= R^{(1)}h_2S(R'{}^{(1)}) \ot S(h_1)R^{(2)}h_3 R'{}^{(2)},
\endalign
$$
for all $h\in H$.
\endproclaim

\demo{Proof} Since, by Proposition~1.3, $R'{}^{(2)}R^{(2)}\ot
S^{(2)}(R^{(1)})S(R'{}^{(1)})= 1\ot 1$, we have
$$
\align
\nu(h)& = R^{(2)}h_2R'{}^{(2)}\ot S(h_1)S(R^{(1)})h_3R'{}^{(1)}\\
& = R^{(2)}h_2R'{}^{(2)}{R'''}^{(2)}{R''}^{(2)}\ot S(h_1)
S(R^{(1)})h_3R'{}^{(1)} S^2({R''}^{(1)})S({R'''}^{(1)})\\
& = R^{(2)}h_2R'{}^{(2)}{R'''}^{(2)}{R''}^{(2)}\ot
S^2({R''}^{(1)}) S(h_1)S(R^{(1)}) h_3R'{}^{(1)}S({R'''}^{(1)})\\
& = R^{(2)}h_2{R''}^{(2)}\ot S^2({R''}^{(1)})S(h_1)S(R^{(1)})h_3.
\endalign
$$
So, the first equality is true. The second one follows in a
similar way.\qed
\enddemo

Let $H$ be a braided Hopf algebra. Since $S(\Z(H))= \Z(H)$, it is
true that
$$
(S\ot S)\xcirc \nu\xcirc S^{-1}(h)\in H\ot \Z(H)\quad \text{for
all $h\in H$}.
$$
Using the expressions for $\nu$ given in Definition~1.1 and
Proposition~1.9 to compute this map, we obtain that for all $h\in
H$:
$$
\align (S\ot S)\xcirc \nu\xcirc S^{-1}(h) &=
R^{(2)}h_2R'{}^{(2)}\ot R^{(1)}h_1
S(R'{}^{(1)})S(h_3)\\
&= R^{(1)}h_2R'{}^{(1)}\ot h_1 R'{}^{(2)}S(h_3)S(R^{(2)})\\
&=R^{(2)}h_2R'{}^{(2)}\ot h_1 S(R'{}^{(1)})S(h_3)S^2(R^{(1)})\\
&=S(R^{(1)})h_2R'{}^{(1)}\ot R^{(2)} h_1 R'{}^{(2)} S(h_3).
\endalign
$$

\proclaim{Proposition 1.10} If $R$ is a semiquasitriangular
structure for a Hopf algebra $H$, then so is $\tau(R^{-1})$.
\endproclaim

\demo{Proof} Conditions~(1)--(4) of Definition~1.1 follow from
standard arguments for quasitriangular Hopf algebras \cite{M2,
Exercise~2.1.3}. Hence, we only check conditions~(5) and (6). By
Proposition~1.3 and the fact that $(H,R)$ satisfies
Proposition~1.9,
$$
\allowdisplaybreaks \align \nu_{(H,\tau(R^{-1}))}(h) &=
S(R^{(1)})h_2S(R'{}^{(1)})\ot S(h_1)S(R^{(2)}) h_3R'{}^{(2)}\\
&= R^{(1)}h_2S(R'{}^{(1)})\ot S(h_1) R^{(2)} h_3R'{}^{(2)}\\
& = \nu_{(H,R)}(h).
\endalign
$$
Consequently, $(H,\tau(R^{-1}))$ satisfies condition~(5).
Moreover, by By Proposition~1.3, Proposition~1.9 and condition~(5)
of Definition~1.1 for $R$, we have
$$
\allowdisplaybreaks
\align
 \nu_{(H,\tau(R^{-1}))}(h)&=\nu_{(H,R)}(h)\\
& =R^{(2)}h_2R'{}^{(2)}\ot S^2(R'{}^{(1)})S(h_1)S(R^{(1)})
h_3\\
& =R^{-(2)}h_2R'{}^{-(2)}\ot S(R'{}^{-(1)})S(h_1)R^{-(1)}h_3.\\
\endalign
$$
So, $(H,\tau(R^{-1}))$ also satisfies condition~(6).\qed
\enddemo

\proclaim{Proposition 1.11} If $R$ is a semiquasitriangular
structure for a Hopf algebra $H$, then $\tau(R)$ and $R^{-1}$ are
semiquasitriangular structures for $H^{\op}$ and $H^{\cop}$.
\endproclaim

\demo{Proof} By Proposition~1.10 it suffices to prove it for
$\tau(R)$. As before, we only check conditions~(5) and (6) of
Definition~1.1, since conditions~(1)--(4) follow from standard
arguments for quasitriangular Hopf algebras. In the rest of the
proof all the multiplications are in $H$. By Proposition~1.3 and
the fact that $(H,R)$ satisfies condition~(5) of Definition~1.1,
we have
$$
\allowdisplaybreaks \align (S^{-1} & \ot H) \xcirc
\nu_{(H^{\op},\tau(R))}(h) = S^{-1}(R'{}^{(1)}h_2R^{(1)})\ot
R'{}^{(2)}h_3S^{-1}(R^{(2)})S^{-1}(h_1)\\
&=S^{-1}(R^{(1)})S^{-1}(h)_2S^{-1}(R'{}^{(1)})\ot
R'{}^{(2)}S(S^{-1}(h)_1)S^{-1}(R^{(2)})S^{-1}(h)_3\\
&=R^{(1)}S^{-1}(h)_2R'{}^{(1)}\ot S(R'{}^{(2)})S(S^{-1}(h)_1)R^{(2)}S^{-1}(h)_3\\
&= \nu_{(H,R)}(S^{-1}(h)),
\endalign
$$
for all $h\in H$. From this it follows immediately that
$\nu_{(H^{\op},\tau(R))}(h) \in H^{\op} \ot \Z(H^{\op})$. Now, let
$h\in H$. Using Propositions~1.3 and 1.9, we obtain
$$
\allowdisplaybreaks \align
\nu_{(H,R)}(S^{-1}(h))&=R^{(2)}S^{-1}(h)_2R'{}^{(2)}  \ot
S(S^{-1}(h)_1)S(R^{(1)})S^{-1}(h)_3R'{}^{(1)}\\
&=S^{-1}(R^{(2)})S^{-1}(h_2)  S^{-1}(R'{}^{(2)})\ot
h_3R^{(1)}S^{-1}(h_1)S^{-1}(R'{}^{(1)})\\
&=S^{-1}(R'{}^{(2)}h_2R^{(2)}) \ot
h_3R^{(1)}S^{-1}(h_1)S^{-1}(R'{}^{(1)}).
\endalign
$$
From these facts it follows easily that $(H^{\op},\tau(R))$
satisfies condition~(6) of Definition~1.1. It remains to prove
that $\tau(R)$ is a semiquasitriangular structure for $H^{\cop}$.
We leave this task to the reader.\qed
\enddemo

\head  2. Normal modules \endhead
Let $(H,R)$ be a semiquasitriangular Hopf algebra.In this section
we introduce the category ${}_{(H,R)}\Cal M^{(H,R)}(n)$ of
left-right normal $(H,R)$-modules and begin the study of its
properties. We assert that this is the suitable category of
representations of $(H,R)$. Evidence that this assertion is
``right'' is given by the facts (proved in this section) that

\smallskip

\itemitem{-} $H$ is a left-right normal $(H,R)$-module in a natural
sense,

\smallskip

\itemitem{-} ${}_{(H,R)}\Cal M^{(H,R)}(n)$ is a braided category,

\smallskip

\itemitem{-}  The full braided subcategory of ${}_{(H,R)}\Cal
M^{(H,R)}(n)$ made out of the finite dimensional modules is rigid,

\smallskip

\itemitem{-} if $(H,R)$ is quasitriangular, then the category of
left $H$-modules is, in a natural way, a braided subcategory of
${}_{(H,R)}\Cal M^{(H,R)}(n)$.

\smallskip

For an exposition of the theory of braided Hopf algebras, we remit
to \cite{J-S}, \cite{Ka}, \cite{Ch-P} and \cite{M2}.

\definition{Definition 2.1} Let $(H,R)$ be an semiquasitriangular
Hopf algebra.

\smallskip

\ni A left-right normal $(H,R)$-module is a vector space $M$,
endowed with a left action $\rho_M\:H\ot M\to M$ and a right
coaction $\nu_M\:M\to M\ot H$, such that for all $h\in H$ and
$m\in M$,

\smallskip

\roster

\item $\nu_M(m)\in M\ot \Z(H)$,

\smallskip

\item $\nu_M(h\cdot m)= R^{(2)}h_2R'{}^{(2)}\cdot m_0\ot S(h_1) S(R^{(1)})
h_3R'{}^{(1)} m_1$.

\endroster

\smallskip

\ni A left normal $(H,R)$-module is a vector space $M$, endowed
with a left action $\rho_M\:H\ot M\to M$ and a left coaction
$\nu_M\:M\to H\ot M$, such that for all $h\in H$ and $m\in M$,

\smallskip

\roster

\item $\nu_M(m)\in \Z(H)\ot M $,

\smallskip

\item $\nu_M(h\cdot m)= S^{-1}(h_3) S^{-1}(R^{(2)})h_1R'{}^{(2)} m_{-1}
\ot R^{(1)}h_2R'{}^{(1)}\cdot m_0$.

\endroster

\smallskip

\ni A right normal $(H,R)$-module is a vector space $M$, endowed
with a right action $\rho_M\:M\ot H\to M$ and a right coaction
$\nu_M\:M\to M\ot H$, such that for all $h\in H$ and $m\in M$,

\smallskip

\roster

\item $\nu_M(m)\in M\ot\Z(H)$,

\smallskip

\item $\nu_M(m\cdot h)= m_0\cdot R'{}^{(1)}h_2 R^{(1)} \ot
m_1R'{}^{(2)} h_3 S^{-1}(R^{(2)})S^{-1}(h_1)$.

\endroster

\smallskip

\ni A right-left normal  $(H,R)$-module is a vector space $M$,
endowed with a right action $\rho_M\:M\ot H\to M$ and a left
coaction $\nu_M\:M\to H\ot M$, such that for all $h\in H$ and
$m\in M$,

\smallskip

\roster

\item $\nu_M(m)\in \Z(H)\ot M$,

\smallskip

\item $\nu_M(m\cdot h)= m_{-1}R'{}^{(1)} h_1 S(R^{(1)})S(h_3)
 \ot m_0\cdot R'{}^{(2)}h_2 R^{(2)}$.

\endroster

\smallskip

\noindent In all the cases a morphism $f\:M\to N$ of normal
$H$-modules is a map which is a morphism of $H$-modules and
$H$-comodules.

\enddefinition

For instance, if $H$ is a semicocommutative Hopf algebra, then
each normal $H$-module is a  left-right normal  $(H,1_{H\ot
H})$-module, and if $(H,R)$ is a quasitriangular Hopf algebra,
then each left $H$-module can be think as a left-right normal
$(H,R)$-module with trivial coaction.

\smallskip

From now on by a normal $(H,R)$-module we understand a left-right
normal $(H,R)$-module.

\proclaim{Lemma 2.2}  Let $(H,R)$ be an semiquasitriangular Hopf
algebra. For each $h\in H$,
$$
(H\ot \Delta)\xcirc\nu(h)\in H\ot \Z(H)\ot H.
$$
\endproclaim

\demo{Proof} By item~(1) of Definition~1.1, we have:
$$
\allowdisplaybreaks \align (H \ot\Delta)& \xcirc\nu(h) =
(H\ot\Delta)\bigl(R^{(2)}h_2R'{}^{(2)}
\ot S(h_1) S(R^{(1)})h_3R'{}^{(1)}\bigr)\\
&=R^{(2)}h_3R'{}^{(2)}\ot S(h_2)S(R^{(1)}_2)h_4 R'{}^{(1)}_1\ot
S(h_1) S(R^{(1)}_1)h_5R'{}^{(1)}_2\\
&=R^{(2)}\wt{R}{}^{(2)}h_3R'{}^{(2)}{\wt{R}'}{}^{(2)}\ot S(h_2)
S(\wt{R}{}^{(1)})h_4R'{}^{(1)}\ot S(h_1)S(R^{(1)})h_5
{\wt{R}'}{}^{(1)}.
\endalign
$$
Now, the assertion follows immediately from item~(5) of
Definition~1.1.\qed
\enddemo

\proclaim{Proposition 2.3}  If $(H,R)$ is an semiquasitriangular
Hopf algebra, then $H$, endowed with the left regular action and
the coaction $\nu$ introduced in Definition~1.1, is a normal
$H$-module.
\endproclaim

\demo{Proof} It is immediate that $\nu$ is a counitary map that
satisfies condition~(1) of Definition~2.1. Next, we check that it
is coassociative. By item~(2) of Proposition~1.2, we have:
$$
\allowdisplaybreaks \align (\nu\ot H)\xcirc \nu(h)& =(\nu\ot
H)\bigl(R^{(2)}h_2\ov{R}^{(2)}
\ot S(h_1) S(R^{(1)})h_3\ov{R}^{(1)}\bigr)\\
& =\wt{R}{}^{(2)}R^{(2)}_2h_3\ov{R}^{(2)}_2\wh{R}{}^{(2)}\ot
S(R^{(2)}_1 h_2\ov{R}^{(2)}_1) S(\wt{R}{}^{(1)})
R^{(2)}_3 h_4 \ov{R}^{(2)}_3\wh{R}{}^{(1)}\\
&\phantom{= \,\,} \ot S(h_1) S(R^{(1)})h_5\ov{R}^{(1)}\\
&=\wt{R}{}^{(2)}R'{}^{(2)}h_3{\ov{R}'}^{(2)}\wh{R}{}^{(2)}\ot
S({R''}^{(2)}h_2{\ov{R}''}^{(2)})S(\wt{R}{}^{(1)}){R}^{(2)}
h_4\ov{R}^{(2)}\wh{R}{}^{(1)}\\
&\phantom{= \,\,} \ot S(h_1)S(R^{(1)}R'{}^{(1)}{R''}^{(1)})h_5
\ov{R}^{(1)}{\ov{R}'}^{(1)}{\ov{R}''}^{(1)}.
\endalign
$$
Hence, by conditions~(1) and (4) of Definition~1.1 and
Proposition~1.3, we have:
$$
\allowdisplaybreaks \align (\nu\ot H)\xcirc
\nu(h)&=\wt{R}{}^{(2)}h_3{\ov{R}'}^{(2)}\ot S({R''}^{(2)}
h_2{\ov{R}''}^{(2)})S(\wt{R}{}^{(1)}_1){R}^{(2)}h_4
\ov{R}^{(2)}{\ov{R}'}^{(1)}_2\\
&\phantom{= \,\,} \ot S(h_1)S(R^{(1)}\wt{R}{}^{(1)}_2
{R''}^{(1)})h_5\ov{R}^{(1)}{\ov{R}'}^{(1)}_1{\ov{R}''}^{(1)}\\
&=\wt{R}{}^{(2)}h_3{\ov{R}'}^{(2)}\ot
S(\wt{R}{}^{(1)}_1{R''}^{(2)}h_2
{\ov{R}''}^{(2)}){R}^{(2)}h_4\ov{R}^{(2)}{\ov{R}'}^{(1)}_2\\
&\phantom{= \,\,} \ot S(R^{(1)}\wt{R}{}^{(1)}_2{R''}^{(1)}h_1)h_5
\ov{R}^{(1)}{\ov{R}'}^{(1)}_1{\ov{R}''}^{(1)}\\
&= \wt{R}{}^{(2)}h_3{\ov{R}'}^{(2)}\ot
S({R''}^{(2)}\wt{R}{}^{(1)}_2
h_2{\ov{R}''}^{(2)})R^{(2)}h_4{\ov{R}'}^{(1)}_1\ov{R}^{(2)}\\
&\phantom{= \,\,} \ot S(R^{(1)}{R''}^{(1)}\wt{R}{}^{(1)}_1h_1)h_5
{\ov{R}'}^{(1)}_2\ov{R}^{(1)}{\ov{R}''}^{(1)}\\
&= \wt{R}{}^{(2)}h_3{\ov{R}'}^{(2)}\ot S({\ov{R}''}^{(2)})S(h_2)
S(\wt{R}{}^{(1)}_2)S({R''}^{(2)})R^{(2)}h_4{\ov{R}'}^{(1)}_1
\ov{R}^{(2)}\\
&\phantom{= \,\,} \ot S(h_1)S(\wt{R}{}^{(1)}_1)S({R''}^{(1)})
S(R^{(1)})h_5 {\ov{R}'}^{(1)}_2\ov{R}^{(1)}{\ov{R}''}^{(1)}\\
&= \wt{R}{}^{(2)}h_3{\ov{R}'}^{(2)}\ot S({\ov{R}''}^{(2)})S(h_2)
S(\wt{R}{}^{(1)}_2) h_4{\ov{R}'}^{(1)}_1\ov{R}^{(2)}\\
&\phantom{= \,\,} \ot S(h_1)S(\wt{R}{}^{(1)}_1)h_5
{\ov{R}'}^{(1)}_2 \ov{R}^{(1)}{\ov{R}''}^{(1)}.
\endalign
$$
On the other hand, by Proposition~1.3 and Lemma~2.2,
$$
\allowdisplaybreaks \align (H\ot \Delta)\xcirc \nu(h) & =
\wt{R}{}^{(2)}h_3{\ov{R}'}^{(2)}\ot S(h_2)
S(\wt{R}{}^{(1)}_2)h_4{\ov{R}'}^{(1)}_1\ot
S(h_1)S(\wt{R}{}^{(1)}_1)h_5
{\ov{R}'}^{(1)}_2\\
&= \wt{R}{}^{(2)}h_3{\ov{R}'}^{(2)}\ot S(h_2)S(\wt{R}{}^{(1)}_2)
h_4{\ov{R}'}^{(1)}_1S(\ov{R}^{(2)}{\ov{R}''}^{(2)})\\
&\phantom{= \,\,} \ot S(h_1)S(\wt{R}{}^{(1)}_1) h_5
{\ov{R}'}^{(1)}_2 S(\ov{R}^{(1)}){\ov{R}''}^{(1)}\\
&= \wt{R}{}^{(2)}h_3{\ov{R}'}^{(2)}\ot S(h_2)S(\wt{R}{}^{(1)}_2)
h_4{\ov{R}'}^{(1)}_1S({\ov{R}''}^{(2)})\ov{R}^{(2)}\\
&\phantom{= \,\,} \ot S(h_1)S(\wt{R}{}^{(1)}_1) h_5
{\ov{R}'}^{(1)}_2 \ov{R}^{(1)}{\ov{R}''}^{(1)}\\
&= \wt{R}{}^{(2)}h_3{\ov{R}'}^{(2)}\ot S({\ov{R}''}^{(2)})S(h_2)
S(\wt{R}{}^{(1)}_2) h_4{\ov{R}'}^{(1)}_1\ov{R}^{(2)}\\
&\phantom{= \,\,} \ot S(h_1)S(\wt{R}{}^{(1)}_1)h_5
{\ov{R}'}^{(1)}_2 \ov{R}^{(1)}{\ov{R}''}^{(1)}.
\endalign
$$
To finish the proof only remain to check that $\nu$ satisfies
condition~(2) of Definition~2.1. But, by Proposition~1.3 and
item~(5) of Definition~1.1,
$$
\allowdisplaybreaks \align \nu(hl)&=
R^{(2)}h_2l_2{\ov{R}'}^{(2)}\ot S(l_1)S(h_1)
S(R^{(1)}) h_3l_3{\ov{R}'}^{(1)}\\
&=R^{(2)}h_2R'{}^{(2)} \ov{R}^{(2)}l_2{\ov{R}'}^{(2)}\ot S(l_1)
S(h_1) S(R^{(1)})h_3R'{}^{(1)} S(\ov{R}^{(1)})l_3{\ov{R}'}^{(1)}\\
&=R^{(2)}h_2R'{}^{(2)} \ov{R}^{(2)}l_2{\ov{R}'}^{(2)}\ot S(h_1)
S(R^{(1)})h_3R'{}^{(1)}S(l_1) S(\ov{R}^{(1)})l_3{\ov{R}'}^{(1)}\\
&=\nu(h)\nu(l),
\endalign
$$
as desired.\qed
\enddemo

\proclaim{Corollary 2.4} If $(H,R)$ be a semiquasitriangular Hopf
algebra, then $H$ is a left normal $(H, R)$-module via the left
regular action and the left coaction
$$
\nu_1(h):= S^{-1}(h_3) S^{-1}(R^{(2)})h_1R'{}^{(2)} \ot
R^{(1)}h_2R'{}^{(1)},
$$
it is a right normal $(H,R)$-module via the right regular action
and the right coaction
$$
\nu_2(h):=R'{}^{(1)}h_2 R^{(1)} \ot R'{}^{(2)} h_3
S^{-1}(R^{(2)})S^{-1}(h_1),
$$
and it is a right-left $(H, R)$-module via the right regular
action and the left coaction
$$
\nu_3(h):= R'{}^{(1)} h_1 S(R^{(1)})S(h_3) \ot R'{}^{(2)}h_2
R^{(2)}.
$$
\endproclaim

\demo{Proof} By Proposition~2.3, it suffices to note that a left
normal $(H,R)$-module is the same that a normal
$(H^{\cop},\tau(R))$-module, a right normal $(H,R)$-module is the
same that a normal $(H^{\op},\tau(R))$-module and a right-left
$(H, R)$-module is the same that a normal $(H^{\op
\,\cop},R)$-module.\qed
\enddemo

\proclaim{Theorem 2.5} Let $(H,R)$ be a semiquasitriangular Hopf
algebra. The category ${}_{(H,R)}\Cal M^{(H,R)}(n)$, of normal
$(H,R)$-modules, is a braided category. The unit object is $k$,
endowed with the trivial action and the trivial coaction, and the
tensor product is the usual tensor product over $k$, endowed with
the diagonal action and the diagonal coaction. The associative and
unit constraints are the usual ones and the braid $c$ is given by
$c_{MN}(m\ot n)= R^{(2)}\cdot n_0\ot R^{(1)}n_1\cdot m$.
\endproclaim

\demo{Proof} Let $M$ and $N$ be normal $H$-modules. It is obvious
that $m_0\ot n_0 \ot m_1n_1\in M\ot N\ot \Z(H)$
for all $h\in H$, $n\in N$ and $m\in M$. Moreover, by item~(2) of
Definition~1.11 and items~(2) and (5) of Definition~1.1, we have
$$
\allowdisplaybreaks
\align
\nu(h\cdot (m\ot n)) & =(h_1\cdot m)_0 \ot (h_2\cdot n)_0\ot
(h_1\cdot m)_1(h_2\cdot n)_1\\
& =R^{(2)}h_2R'{}^{(2)}\cdot m_0\ot \ov{R}^{(2)}h_5{\ov{R}'}^{(2)}
\cdot n_0\\
&\phantom{=\,\,} \ot S(h_1)S(R^{(1)})h_3R'{}^{(1)}m_1S(h_4)
S(\ov{R}^{(1)}) h_6{\ov{R}'}^{(1)}n_1\\
& =R^{(2)}h_2R'{}^{(2)}\cdot m_0\ot \ov{R}^{(2)}h_5{\ov{R}'}^{(2)}
\cdot n_0\\
&\phantom{=\,\,}\ot S(h_1)S(R^{(1)})h_3S(h_4)S(\ov{R}^{(1)})h_6
{\ov{R}'}^{(1)}R'{}^{(1)}m_1n_1\\
& =R^{(2)}h_2R'{}^{(2)}\cdot m_0\ot \ov{R}^{(2)}h_3{\ov{R}'}^{(2)}
\cdot n_0\\
&\phantom{=\,\,}\ot S(h_1)S(\ov{R}^{(1)}R^{(1)})h_4{\ov{R}'}^{(1)}
R'{}^{(1)}m_1n_1\\
& =R^{(2)}_1h_2R'{}^{(2)}_1\cdot m_0\ot R^{(2)}_2h_3 R'{}^{(2)}_2
\cdot n_0 \ot S(h_1)S(R^{(1)}) h_4R'{}^{(1)}m_1n_1\\
& =R^{(2)}h_2R'{}^{(2)}\cdot (m_0\ot n_0)\ot S(h_1)S(R^{(1)})
h_3R'{}^{(1)}m_1n_1,
\endalign
$$
for each $h\in H$. Hence, the tensor product $M\ot N$ is a normal
$H$-module via the diagonal action and the diagonal coaction.
Moreover, it is immediate that $k$ is a normal $H$-module, and it
is clear that the usual associative and unit constraints are
$H$-linear and $H$-colinear maps. So, ${}_{(H,R)}\Cal
M^{(H,R)}(n)$ is a monoidal category. To prove that it is a
braided category with braid $c$, we must show that $c$ is a
natural isomorphism of normal $(H,R)$-modules, and that
$$
c_{M\ot N,P}=(c_{MP}\ot N)\xcirc (M\ot c_{NP})\quad\text{and}\quad
c_{M,N\ot P}=(N\ot c_{MP})\xcirc (c_{MN}\ot P)
$$
for all $M,N,P\in {}_{(H,R)}\Cal M^{(H,R)}(n)$. We do this in
several steps.

\medskip

\ni {\it $c_{MN}$ is $H$-linear}:\enspace For $h\in H$, $m\in M$
and $n\in N$, we have:
$$
\allowdisplaybreaks \align
c_{MN}(h\cdot &(m\ot n))=c_{MN}(h_1\cdot m\ot h_2\cdot n)\\
&= R^{(2)}\cdot (h_2\cdot n)_0 \ot R^{(1)} (h_2\cdot n)_1 \cdot (h_1\cdot m)\\
&=R^{(2)}\ov{R}^{(2)}h_3{\ov{R}'}^{(2)}\cdot n_0 \ot
R^{(1)}S(h_2)S(\ov{R}^{(1)})h_4{\ov{R}'}^{(1)}n_1h_1\cdot m\\
&=R^{(2)}\ov{R}^{(2)}h_3{\ov{R}'}^{(2)}\cdot n_0 \ot
R^{(1)}h_1S(h_2)S(\ov{R}^{(1)})h_4{\ov{R}'}^{(1)}n_1\cdot m\\
&=R^{(2)}\ov{R}^{(2)}h_1{\ov{R}'}^{(2)}\cdot n_0 \ot
R^{(1)}S(\ov{R}^{(1)})h_2{\ov{R}'}^{(1)}n_1\cdot m\\
&=h_1{\ov{R}'}^{(2)}\cdot n_0 \ot h_2{\ov{R}'}^{(1)}n_1\cdot m\\
&=h\cdot({\ov{R}'}^{(2)}\cdot n_0 \ot {\ov{R}'}^{(1)}n_1\cdot m)\\
&=h\cdot c_{MN}(m\ot n),
\endalign
$$
where the third equality follows from item~(2) of Definition~2.1,
the fourth one follows from item~(1) of Definition~2.1 and
item~(5) of Definition~1.1 and the sixth one follows from
Proposition~1.3.

\medskip

\ni {\it $c_{MN}$ is $H$-colinear}:\enspace For $m\in M$ and $n\in
N$, we have:
$$
\allowdisplaybreaks \align
\nu&(c_{MN} (m\ot n))= \nu(R^{(2)}\cdot n_0 \ot R^{(1)}n_1\cdot m)\\
&= (R^{(2)}\cdot n_0)_0 \ot (R^{(1)}n_1\cdot m)_0\ot
(R^{(2)}\cdot n_0)_1(R^{(1)}n_1\cdot m)_1\\
&=\ov{R}^{(2)}R^{(2)}_2{\ov{R}'}^{(2)}\cdot n_0 \ot
\wt{R}{}^{(2)}R^{(1)}_2n_3\wt{R}'{}^{(2)}\cdot m_0\ot\\
& \phantom{=\,\,} S(R^{(2)}_1)S(\ov{R}^{(1)})R^{(2)}_3
{\ov{R}'}^{(1)} n_1S(n_2) S(R^{(1)}_1)S(\wt{R}{}^{(1)})
R^{(1)}_3n_4\wt{R}'{}^{(1)}m_1\\
&=\ov{R}^{(2)}R^{(2)}_2{\ov{R}'}^{(2)}\cdot n_0 \ot
\wt{R}{}^{(2)}R^{(1)}_2n_1\wt{R}'{}^{(2)}\cdot m_0\ot\\
& \phantom{=\,\,} S(\ov{R}^{(1)}R^{(2)}_1)R^{(2)}_3
{\ov{R}'}^{(1)} S(\wt{R}{}^{(1)}R^{(1)}_1)R^{(1)}_3n_2
\wt{R}'{}^{(1)}m_1\\
&=R^{(2)}_1\ov{R}^{(2)}{\ov{R}'}^{(2)}\cdot n_0 \ot
R^{(1)}_1\wt{R}{}^{(2)}n_1\wt{R}'{}^{(2)}\cdot m_0\ot\\
& \phantom{=\,\,} S(\ov{R}^{(1)})S(R^{(2)}_2)R^{(2)}_3
{\ov{R}'}^{(1)} S(\wt{R}{}^{(1)})S(R^{(1)}_2)R^{(1)}_3n_2
\wt{R}'{}^{(1)}m_1\\
&=R^{(2)}\ov{R}^{(2)}{\ov{R}'}^{(2)}\cdot n_0 \ot
R^{(1)}\wt{R}{}^{(2)}\wt{R}'{}^{(2)}n_1\cdot m_0\ot
S(\ov{R}^{(1)}){\ov{R}'}^{(1)} S(\wt{R}{}^{(1)})
\wt{R}'{}^{(1)}m_1n_2\\
&=R^{(2)}\cdot n_0 \ot R^{(1)}n_1\cdot m_0\ot m_1n_2\\
&=(c_{MN}\ot H)\xcirc \nu (m\ot n),
\endalign
$$
where the third equality follows from item~(2) of Definition~2.1,
the fifth one follows from items~(3) and (4) of Definition~1.1 and
the seventh one follows from Proposition~1.3.

\medskip

\ni {\it $c$ is a natural isomorphism}:\enspace Let $f\:M\to M'$
and $g\:N\to N'$ morphisms of normal modules. For each $m\in M$
and $n\in N$, we have:
$$
\allowdisplaybreaks \align
(g\ot f)c_{MN}(m\ot n)&=g(R^{(2)}\cdot n_0)\ot f(R^{(1)} n_1\cdot m)\\
&=R^{(2)}\cdot g(n_0)\ot R^{(1)}n_1\cdot f(m)\\
&=R^{(2)}\cdot g(n)_0\ot R^{(1)}g(n)_1\cdot f(m)\\
&=c_{M'N'}(f(m)\ot g(n)).
\endalign
$$
This shows that $c$ is a natural transformation. In order to prove
that $c_{MN}$ is a bijective map it suffices to note that
$c=l_{\tau(R)} \xcirc \ov{c}$, where $\ov{c}\:M\ot N\to N\ot M$
and $l_{\tau(R)}\:N\ot M\to N\ot M$ are the maps defined by
$$
\ov{c}(m\ot n):= n_0 \ot n_1\cdot m \quad\text{and}\quad
l_{\tau(R)}(n\ot m) := R^{(2)}\cdot n\ot R^{(1)}\cdot m,
$$
which clearly are bijective.

\medskip

\ni {\it It is true that $(c_{MP} \ot N)\xcirc (M\ot c_{NP})=
c_{M\ot N,P}$:\enspace by condition~(2) of Definition~2.1,
conditions~(1) and (3) of Definition~1.1 and Proposition~1.3,
$$
\allowdisplaybreaks \align (c_{MP} &\ot N)\xcirc (M\ot
c_{NP})(m\ot n\ot p)= (c_{MP}
\ot N)(m\ot R^{(2)}\cdot p_0 \ot R^{(1)}p_1\cdot n)\\
&= \ov{R}^{(2)}\cdot (R^{(2)}\cdot p_0)_0\ot
\ov{R}^{(1)}(R^{(2)}\cdot p_0)_1\cdot m \ot R^{(1)}p_1\cdot n\\
&= \ov{R}^{(2)}\wt{R}{}^{(2)} R^{(2)}_2\wt{R}'{}^{(2)}\cdot p_0\ot
\ov{R}^{(1)}S(R^{(2)}_1)S(\wt{R}{}^{(1)})R^{(2)}_3\wt{R}'{}^{(1)}
p_1\cdot m \ot R^{(1)}p_2\cdot n\\
&= \ov{R}^{(2)}R^{(2)}_1\wt{R}{}^{(2)}\wt{R}'{}^{(2)}\cdot p_0\ot
\ov{R}^{(1)}S(\wt{R}{}^{(1)})S(R^{(2)}_2)R^{(2)}_3\wt{R}'{}^{(1)}
 p_1\cdot m\ot R^{(1)}p_2\cdot n\\
&= \ov{R}^{(2)}R^{(2)}\cdot p_0\ot \ov{R}^{(1)}p_1\cdot m \ot
R^{(1)}p_2\cdot n\\
&= R^{(2)}\cdot p_0\ot \ov{R}^{(1)}_1p_1\cdot m
\ot R^{(1)}_2p_2\cdot n\\
&= c_{M\ot N,P}(m\ot n\ot p),
\endalign
$$
where the third equality follows from item~(2) of Definition~2.1,
the fourth one follows from item~(3) of Definition~1.1, the fifth
one follows from Proposition~1.3 and the sixth one follows from
item~(1) of Definition~1.1.

\medskip

\ni {\it It is true that $(N\ot c_{MP})\xcirc (c_{MN} \ot P)=
c_{M, N\ot P}$:\enspace by item~(1) of Definition~2.1 and item~(2)
of Definition~1.1,
$$
\allowdisplaybreaks \align (N\ot c_{MP})\xcirc (c_{MN} \ot P)(m\ot
n\ot p)&=(N\ot c_{MP})
(R^{(2)}\cdot n_0\ot R^{(1)} n_1\cdot m \ot p)\\
&=R^{(2)}\cdot n_0\ot \ov{R}^{(2)}\cdot p_0\ot\ov{R}^{(1)}p_1
R^{(1)} n_1\cdot m \\
&=R^{(2)}\cdot n_0\ot\ov{R}^{(2)}\cdot p_0\ot\ov{R}^{(1)}R^{(1)}
n_1p_1\cdot m \\
&=R^{(2)}_1\cdot n_0\ot R^{(2)}_2\cdot p_0\ot R^{(1)}n_1p_1\cdot m
\\
&=c_{M,N\ot P}(m\ot n\ot p).
\endalign
$$
This finish the proof.\qed
\enddemo

\remark{Remark 2.6} The inverse of the braid $c$ introduced in
Theorem~2.5 is given by $c^{-1}_{NM}(n\ot m) =
S(n_1)S(R^{(1)})\cdot m\ot R^{(2)}\cdot n_0$. In fact, let
$\ov{c}\: M\ot N\to N\ot M$ and $l_{\tau(R)}\:N\ot M\to N\ot M$ be
as in the proof of Theorem~2.5. It is easy to see that $\ov{c}^{-1}
(n\ot m) = S(n_1)\cdot m\ot n_0$ and $l_{\tau(R)}^{-1}(n\ot m) =
R^{(2)}\cdot n\ot S(R^{(1)})\cdot m$. Hence, by item~(5) of
Definition~1.1, item~(3) of Proposition~1.2 and Proposition~1.3, we have
$$
\align
c^{-1}_{NM}(n\ot m) & = \ov{c}^{-1}\bigl(l_{\tau(R)}^{-1}(n\ot
m)\bigr)\\
& = S\bigl((R^{(2)}\cdot n)_1\bigr) S(R^{(1)})\cdot m\ot (R^{(2)}\cdot n)_0\\
&= S\bigl(S(R^{(2)}_1)S(\ov{R}{}^{(1)})R^{(2)}_3\ov{R}'{}^{(1)}
n_1\bigr) S(R^{(1)})\cdot m\ot\ov{R}{}^{(2)}R^{(2)}_2
\ov{R}'{}^{(2)} \cdot n_0\\
&= S\bigl(S(\ov{R}{}^{(1)})S(R^{(2)}_2)R^{(2)}_3\ov{R}'{}^{(1)}
n_1\bigr)S(R^{(1)})\cdot m\ot R^{(2)}_1\ov{R}{}^{(2)}
\ov{R}'{}^{(2)}\cdot n_0\\
&= S\bigl(S(\ov{R}{}^{(1)})\ov{R}'{}^{(1)}n_1\bigr)S(R^{(1)})\cdot m
\ot R^{(2)}\ov{R}{}^{(2)}\ov{R}'{}^{(2)}\cdot n_0\\
&= S(n_1)S(R^{(1)})\cdot m \ot R^{(2)}\cdot n_0,
\endalign
$$
as we assert.
\endremark

\proclaim{Corollary 2.7} Let $(H,R)$ be a semiquasitriangular Hopf
algebra. The categories ${}_{(H,R)}^{(H,R)}\Cal M(n)$, of left
normal $(H,R)$-modules, $\Cal M_{(H,R)}^{(H,R)}(n)$, of right
normal $(H,R)$-modules, and ${}^{(H,R)}\Cal M_{(H,R)}(n)$, of
right-left normal $(H,R)$-modules, are braided categories. In all
the cases the unit object is $k$, endowed with the trivial action
and the trivial coaction; the tensor product is the usual tensor
product over $k$, endowed with the diagonal action and the
diagonal coaction and the associative and unit constraints are the
usual ones. The braids are given respectively by ${c_1}_{MN} (m\ot
n)= R^{(1)}\cdot n_0\ot R^{(2)}n_{-1}\cdot m$, ${c_2}_{MN}
(m\ot n)= n_0 \cdot R^{(1)} \ot m \cdot n_1 R^{(2)}$ and
${c_3}_{MN} (m\ot n)= n_0 \cdot R^{(2)} \ot m \cdot n_1
R^{(1)}$.
\endproclaim

\demo{Proof} Proceed as in the proof of Corollary~2.4.\qed
\enddemo

Recall that an object $V$ of a braided category $\Cal C$ is rigid
if there exists an object $V^*$, endowed with arrows $\ev_V:V^*\ot
V\to 1_{\Cal C}$ and $\coev_V: 1_{\Cal C}\to V\ot V^*$, where
$1_{\Cal C}$ is the unit object of $\Cal C$, satisfying
$$
\ide_V= (V\ot \ev_V)\xcirc (\coev_V\ot V)\quad\text{and}\quad
\ide_{V^*}=(\ev_V\ot V^*)\xcirc (V^*\ot \coev_V).
$$
The object $V^*$, which is unique unless a canonical isomorphism,
is called the left dual of $V$, and the morphisms $\ev_V$ and
$\coev_V$ are called the evaluation and the coevaluation maps of
$V$, respectively. Let $U,V$ rigid objects of $\Cal C$ and let
$f\:U\to V$ be a map of $\Cal C$. The transpose map $f^*\:V^*\to
U^*$ of $f$ is defined by
$$
f^* := (\ev_V\ot U^*)\xcirc (V^*\ot f \ot U^*)\xcirc (V^*\ot
\coev_U).
$$
A braided category is said to be rigid if each object has a left
dual. Let $(H,R)$ be a semiquasitriangular Hopf algebra. We are
going to prove that the category of finite dimensional left-right
normal $(H,R)$- modules is a rigid braided category.

\smallskip

Let $M$ be a finite dimensional left-right normal $(H,R)$-module.
Given $f\in M^*$ and $h\in H$ we define $h\cdot f$ by $(h\cdot
f)(m) = f(S(h)\cdot m)$ and we define $\nu_{M^*}(f)\in M^*\ot
\Z(H)$ by $\nu(f) = \sum_i f({m_i}_0)m_i^* \ot S^{-1}({m_i}_1)$,
where $\{m_i,m_i^*\}_{i\in I}$ are dual bases of $M$.

\proclaim{Theorem 2.8} $M^*$ is a left-right normal
$(H,R)$-module.
\endproclaim

\demo{Proof} It is immediate that $h\ot f\mapsto h\cdot f$ is an
action and that $(M^*\ot \epsilon)\xcirc \nu(f) = f$. Let us see
that $(M^*\ot \De)\xcirc \nu = (\nu\ot H)\xcirc \nu$. By
definition
$$
\align (M^*\ot \De)\xcirc \nu_{M^*}(f) & = \sum_{i\in I}
f({m_i}_0)m_i^* \ot S^{-1}({m_i}_2)\ot S^{-1}({m_i}_1)\\
\intertext{and}
(\nu\ot H)\xcirc \nu_{M^*}(f) & = \sum_{i,j \in I} f({m_i}_0)
m_i^*({m_j}_0) m_j^* \ot S^{-1}({m_j}_1)\ot S^{-1}({m_i}_1).
\endalign
$$
Evaluating in the first factor of these expressions in $m_k$ for
$1\le k\le n$, we reduce to prove that
$$
f({m_k}_0)\ot S^{-1}({m_k}_2)\ot S^{-1}({m_k}_1) = \sum_{i\in I}
f({m_i}_0) m_i^*({m_k}_0)\ot S^{-1}({m_k}_1)\ot S^{-1}({m_i}_1).
$$
This follows to applying $(f\ot S^{-1}\ot S^{-1})\xcirc (N\ot
\tau)\xcirc (\nu_M\ot H)$ to the equality
$$
{m_k}_0\ot {m_k}_1 = \sum_{i\in I} {m_i}_0m_i^*({m_k}_0)\ot
{m_k}_1.
$$
It remains to prove that condition~(2) of Definition~2.1 is
satisfied. We must see that
$$
(h\cdot f)_0(m)h^*((h\cdot f)_1) = (R^{(2)}h_2R'{}^{(2)}\cdot
f_0)(m) h^*\bigl(S(h_1)S(R^{(1)})h_3R'{}^{(1)}f_1\bigr)
$$
for all $m\in M$ and $h^*\in H^*$. On one hand we have
$$
\align (h\cdot f)_0(m)h^*((h\cdot f)_1) & = \sum_{i\in I}
f(S(h){m_i}_0) m_i^*(m)m^*(S^{-1}({m_i}_1))\\
& = \sum_{i\in I} f\bigl(S(h){m_i}_0m_i^*(m)\bigr) h(S^{-1}({m_i}_1))\\
& = f\bigl(S(h)m_0) h^*(S^{-1}(m_1)),
\endalign
$$
where the last equality follows from the fact that $\sum_{i\in I}
m_i m_i^*(m) = m$. On the other hand,
$$
\align
& (R^{(2)}h_2R'{}^{(2)}\cdot f_0)(m)h^*\bigl(S(h_1)S(R^{(1)})h_3
R'{}^{(1)}f_1\bigr)\\
& = \sum_{i\in I} f({m_i}_0)
m_i^*\bigl(S(R^{(2)}h_2R'{}^{(2)})\cdot m\bigr)
h^*\bigl(S(h_1)S(R^{(1)})h_3R'{}^{(1)} S^{-1}({m_i}_1)\bigr)\\
& = \sum_{i\in I} f\bigl({m_i}_0m_i^*\bigl(S(R^{(2)}h_2R'{}^{(2)})\cdot m\bigr)
\bigr)h^*\bigl(S(h_1)S(R^{(1)})h_3R'{}^{(1)} S^{-1}({m_i}_1)\bigr)\\
& = f\bigl(\bigl(S(R^{(2)}h_2R'{}^{(2)})\cdot
m\bigr)_0\bigr)h^*\bigl(S(h_1)S(R^{(1)})
h_3R'{}^{(1)}S^{-1}\bigl(\bigl(S(R^{(2)}h_2R'{}^{(2)})\cdot m\bigr)_1\bigr)\bigr)\\
& = f\bigl(\ov{R}{}^{(2)} S(R^{(2)}_2h_3R'{}^{(2)}_2)\ov{R}'{}^{(2)}\cdot m_0
\bigr)h^*\Bigl(S(h_1)S(R^{(1)})h_5R'{}^{(1)}\\
&\quad\,\, S^{-1}\bigl(S^2(R^{(2)}_3h_4 R'{}^{(2)}_3)
S(\ov{R}{}^{(1)})S(R^{(2)}_1
h_2R'{}^{(2)}_1)\ov{R}'{}^{(1)}  m_1\bigr)\Bigr)\\
& = f\bigl(\ov{R}{}^{(2)}S(R'{}^{(2)}_2)S(h_3)S(R^{(2)}_2)
\ov{R}'{}^{(2)}\cdot m_0 \bigr)h^*\Bigl(S(h_1)S(R^{(1)})h_5R'{}^{(1)}\\
&\quad\,\,S^{-1}(m_1)S^{-1}(\ov{R}'{}^{(1)})R^{(2)}_1
h_2R'{}^{(2)}_1 \ov{R}{}^{(1)}S(R'{}^{(2)}_3)S(h_4)S(R^{(2)}_3)\Bigr)\\
\allowdisplaybreak
& = f\bigl(S(R'{}^{(2)}_2\ov{R}{}^{(2)})S(h_3)S(\ov{R}'{}^{(2)}R^{(2)}_2)
\cdot m_0 \bigr)h^*\Bigl(S(h_1)S(R^{(1)})h_5R'{}^{(1)}\\
&\quad\,\,S^{-1}(m_1)\ov{R}'{}^{(1)}R^{(2)}_1
h_2R'{}^{(2)}_1 S(R'{}^{(2)}_3\ov{R}{}^{(1)})S(h_4)S(R^{(2)}_3)\Bigr)\\
& =
f\bigl(S(\ov{R}{}^{(2)}R'{}^{(2)})S(h_3)S(\ov{R}'{}^{(2)}R^{(2)}_2)
\cdot m_0 \bigr)h^*\Bigl(S(h_1)S(R^{(1)})h_5R'{}^{(1)}\\
&\quad\,\,S^{-1}(m_1)\ov{R}'{}^{(1)}R^{(2)}_1 h_2
S(\ov{R}{}^{(1)})
S(h_4)S(R^{(2)}_3)\Bigr)\\
& =
f\bigl(S(R'{}^{(2)})S(\ov{R}{}^{(2)})S(h_3)S(\ov{R}'{}^{(2)})S(R^{(2)}_1)
\cdot m_0 \bigr)h^*\Bigl(S(h_1)S(R^{(1)})h_5R'{}^{(1)}\\
&\quad\,\,S^{-1}(m_1)R^{(2)}_2\ov{R}'{}^{(1)} h_2
S(\ov{R}{}^{(1)}) S(h_4)S(R^{(2)}_3)\Bigr),
\endalign
$$
where the third equality follows from the fact that $\sum_{i\in I}
m_i m_i^*(m) = m$, the fourth one follows from the fact that
$\nu_M$ satisfies condition~(2) of Definition~2.1, the sixth one
follows from Proposition~1.3 and the seventh and eighth ones
follow from item~(3) of Proposition~1.2. Since, by the discussion
that follows Proposition~1.9,
$$
h_1\ot S(\ov{R}{}^{(2)})S(h_3)S(\ov{R}'{}^{(2)})\ot
\ov{R}'{}^{(1)} h_2 S(\ov{R}{}^{(1)}) S(h_4)\in H\ot H\ot \Z(H),
$$
we obtain
$$
\align
(R^{(2)}h_2R'{}^{(2)}& \cdot f_0)(m)h^*\bigl(S(h_1)S(R^{(1)})h_3
R'{}^{(1)}f_1\bigr)\\
& = f\bigl(S(\ov{R}{}^{(2)}R'{}^{(2)})S(h_3)S(R^{(2)}
\ov{R}'{}^{(2)}) \cdot m_0 \bigr)\\
&\quad\,\,h^*\Bigl(S(h_1)S(R^{(1)})h_5R'{}^{(1)}S^{-1}(m_1)\ov{R}'{}^{(1)}
h_2 S(\ov{R}{}^{(1)})S(h_4)\Bigr)\\
\allowdisplaybreak
& =f\bigl(S(\ov{R}{}^{(2)}R'{}^{(2)})S(h_3)S(R^{(2)}\ov{R}'{}^{(2)})
\cdot m_0 \bigr)\\
&\quad\,\,h^*\Bigl(S(h_1)S(R^{(1)})\ov{R}'{}^{(1)}
h_2 S(\ov{R}{}^{(1)})S(h_4)h_5R'{}^{(1)}S^{-1}(m_1)\Bigr)\\
& =f\bigl(S(\ov{R}{}^{(2)}R'{}^{(2)})S(h_3)\cdot m_0
\bigr)h^*\bigl(S(h_1)
h_2 S(\ov{R}{}^{(1)})R'{}^{(1)}S^{-1}(m_1)\bigr)\\
& =f(S(h)\cdot m_0)h^*(S^{-1}(m_1)),
\endalign
$$
where the third and fourth one equalities follows from
Proposition~1.3, as we need.\qed.

\enddemo

\proclaim{Theorem 2.9} Let $M$ be a finite dimensional left-right
normal $(H,R)$-module. Then, the left-right normal $(H,R)$-module
$M^*$, together with the usual evaluation and coevaluation maps
$\ev_M\:M^*\ot M\to k$ and $\coev_M\:k \to M\ot M^*$, is a left
dual of $M$.
\endproclaim

\demo{Proof} It is well known that $\ev_M$ and $\coev_M$ are
$H$-linear maps and that
$$
\ide_M= (M\ot \ev_M)\xcirc (\coev_M\ot M)\quad\text{and}\quad
\ide_{M^*}=(\ev_M\ot M^*)\xcirc (M^*\ot \coev_M).
$$
It remains to prove that  $\ev_M$ and $\coev_M$ are $H$-colinear
maps. Let $\{m_i,m_i^*\}_{i\in I}$ be dual bases of $M$ and let
$m\in M$. Since $\sum_{i\in I} m_im^*(m_0) \ot m_1=m_0\ot m_1$ and
$m_0\ot m_1\in H\ot \Z(H)$, we have
$$
\align (\ev_M\ot H)\xcirc \nu(f) & = (\ev_M\ot H)(f({m_i}_0)m_i^*
\ot m_0 \ot m_1S^{-1}({m_i}_1))\\
& = f({m_i}_0) m_i^*(m_0) \ot m_1S^{-1}({m_i}_1)\\
&  = f(m_0)\ot m_2S^{-1}(m_1)\\
& = f(m)\ot 1 \endalign
$$
for each $f\in M^*$. So, $\ev_M$ is a morphism of comodules. To
check that $\coev_M$ is also, it suffices to note that since
$$
\sum_{i\in I} {m_i}_0\ot {m_i}_1 \ot m_i^* =\sum_{ij\in I} m_j \ot
{m_i}_1 \ot m_j^*({m_i}_0) m_i^*,
$$
we have
$$
\align
\nu\xcirc \coev(1) &= \sum_{j\in I} \nu(m_j\ot m_j^*)\\
&= \sum_{ij\in I} {m_j}_0 \ot m_j^*({m_i}_0) m_i^* \ot S^{-1}({m_i}_1){m_j}_1\\
&= \sum_{i\in I} {m_i}_0 \ot m_i^* \ot S^{-1}({m_i})_2 {m_i}_1\\
&=\sum_{i\in I} m_i\ot m_i^*\ot 1\\
&=(\coev\ot H)\xcirc \nu(1),
\endalign
$$
as desired.\qed
\enddemo

\proclaim{Corollary 2.10} The category ${}_{(H,R)}\Cal
M_{\fin}^{(H,R)}(n)$, of left-right finite dimensional normal
$(H,R)$-modules, is rigid.
\endproclaim

Let $(H,R)$ be a finite dimensional semiquasitriangular Hopf
algebra and let $\{h_i,h_i^*\}_{i\in I}$ be dual bases of $H$ and
$H^*$. Let $H\bowtie H^*$ denote the tensor product $H\ot H^*$,
endowed with the multiplication $(h\bowtie \psi)(l\bowtie \phi) =
hR^{(2)}l_2R'{}^{(2)} \bowtie \bigl(\psi \leftharpoonup
S(l_1)S(R^{(1)})l_3R'{}^{(1)}\bigr)\phi$ and the codiagonal
comultiplication. We write $h\bowtie\psi$ to denote the element
$h\ot\psi$ of $H\bowtie H^*$. Let $T\in (H\bowtie H^*)\ot
(H\bowtie H^*)$ be the element $T := \sum_{i\in I} (R^{(1)}h_i
\bowtie \ep)\ot (R^{(2)}\bowtie h_i^*)$.

\proclaim{Theorem 2.11} $H\bowtie H^*$ is a Hopf algebra with unit
$1\bowtie \ep$, counit $\ep_H\ot \ep_{H^*}$ and antipode
$S_{H\bowtie H^*}(h\bowtie\psi):= (1\bowtie \phi\xcirc
S)(S(h)\bowtie \ep)$. Moreover the category of left
representations of $H\bowtie H^*$ coincide with the category of
normal $(H,R)$-modules and $(H\bowtie H^*,T)$ is a quasitriangular
Hopf Algebra.

\endproclaim

\demo{Proof} This result can be proved using \cite{Ch-P,
Theorem~5.1.11}, but here we prefer to give a direct proof. Let
$\chi\:H^*\ot H\to H\ot H^*$ be the map given by $\chi(\psi\ot h)
=  (1\bowtie h)(\psi\bowtie 1)$. In order to prove that $H\bowtie
H^*$ is an associative algebra with unit $1\bowtie \ep$, it
suffices to check that $\chi$ is a twisted map in the sense of
\cite{C-V-S}. That is,
$$
\align & \chi\xcirc (H^*\ot\mu) = (\mu\ot H^*)\xcirc(H\ot\chi)
\xcirc (\chi\ot H),\\
& \chi\xcirc (\mu_{H^*}\ot H) = (H\ot \mu_{H^*})\xcirc
(\chi\ot H^*)\xcirc (H^*\ot\chi),\\
&\chi(\ep \ot h) = h\ot\ep\qquad\text{and}\qquad \chi(\psi \ot 1)
= 1\ot\psi.
\endalign
$$
We leave this to the reader. Now, it is immediate that the
category of left representations of $H\bowtie H^*$ coincide with
the category of normal $(H,R)$-modules. In fact, if $M$ is a
normal $(H,R)$-module, then
$$
(h\bowtie\psi)\cdot m = h\cdot m_0\psi(m_1) \quad\text{for $m\in
M$.}
$$
Since the category of finite dimensional normal $(H,R)$-modules is
rigid monoidal, it is true that $H\bowtie H^*$ is a Hopf algebra
with comultiplication  and antipode given by
$$
\De_{H\bowtie H^*}(h\bowtie\psi) = (h\bowtie\psi)\cdot
[(1\bowtie\ep)\ot (1\bowtie\ep)] = (h_1\bowtie \psi_1)\ot
(h_2\bowtie \psi_2)
$$
and
$$
S_{H\bowtie H^*}(h\bowtie\psi) = (1\bowtie S(\phi))(S(h)\bowtie
\ep).
$$
Finally, since the category of normal $(H,R)$-modules is braided,
$$
\align
\tau\xcirc c((1\bowtie\ep)\ot (1\bowtie\ep)) & = \tau\left(
\sum_i R^{(2)}\cdot (1\bowtie h_i^*)\ot R^{(1)}h_i\cdot
(1\bowtie \ep)\right)\\
& = \sum_i (R^{(1)}h_i\bowtie \ep)\ot (R^{(2)}\bowtie h_i^*),
\endalign
$$
is an $R$-structure of $H\bowtie H^*$.\qed

\enddemo

\example{Example 2.12} When $(H,R)$ is quasitriangular, then
$H\bowtie H^* = H\ot H^*$.
\endexample

\example{Example 2.13} Let $G$ be a finite group, $H = k[G]^*$ and
$R = 1\ot 1$. For each $x\in G$ let $\de_x\:G\to k$ be the map
$\de_x(y) = \de_{x,y}$, where $\de_{x,y}$ is the Kronecker symbol.
Then $H\bowtie H^*$ is the tensor product of $k[G]^*$ with $k[G]$,
endowed with the multiplication given by
$$
(\de_x\bowtie y)(\de_{x'}\bowtie y') = \cases \de_x\bowtie yy'
& \text{ if $x' = yxy^{-1}$,}\\ 0 &\text{ in other case.}\endcases
$$
In this case, the $R$-matrix $T$ is $\sum_{x,y\in G}
(\de_x\bowtie 1)\ot (\de_y\bowtie x)$.
\endexample

\head 3. The Drinfeld element of a semiquasitriangular Hopf
algebra \endhead
In this section we show that the properties of the Drinfeld
element of a quasitriangular Hopf algebra remain valid in the
semiquasitriangular setting. However, in this last case same
formulas are more involved  (see for instance Proposition~3.2).

\definition{Definition 3.1} Let $(H, R)$ be a semiquasitriangular
Hopf algebra. The Drinfeld element of $(H,R)$ is the element
$u:=S(R^{(2)})R^{(1)}$ of $H$.
\enddefinition

\proclaim{Proposition 3.2} Assume that $(H, R)$ is a
semiquasitriangular Hopf algebra. Let $T\:H\to H$ be the map
defined by
$$
T(h):=R^{(2)}h_2R'{}^{(2)} S(h_1) S(R^{(1)}) h_3R'{}^{(1)}.
$$
The Drinfeld element $u$ is invertible with inverse
$R^{(2)}S^2(R^{(1)})$. Moreover, $S^2(h)= uT(h)u^{-1}$ for all
$h\in H$.

\endproclaim

\demo{Proof} By conditions~(5) and (6) of Definition~1.1 and
Proposition~1.3, we have
$$
\allowdisplaybreaks \align
S(h_2)uT(h_1)&=S(h_4)uR^{(1)}h_2R'{}^{(1)}
S(R'{}^{(2)})S(h_1)R^{(2)}h_3\\
&=S(R'{}^{(2)})S(h_1)R^{(2)}uR^{(1)}h_2R'{}^{(1)}\\
&=S(R'{}^{(2)})S(h_1)h_2R'{}^{(1)}\\
&=u\epsilon(h).
\endalign
$$
Hence, $S^2(h)u=S^2(h_2)\epsilon(h_1)u=S^2(h_3)
S(h_2)uT(h_1)=uT(h)$. It remain to check that $u$ is invertible
and $u^{-1}=R^{(2)}S^{2}(R^{(1)})$. By the formula proved above,
condition~(4) of Definition~1.1 and Proposition~1.3, we have
$$
\allowdisplaybreaks \align \wt{R}{}^{(2)}S^2(\wt{R}{}^{(1)})u
&=\wt{R}{}^{(2)}uR^{(2)}\wt{R}{}^{(1)}_2R'{}^{(2)}
S(\wt{R}{}^{(1)}_1)S(R^{(1)})\wt{R}{}^{(1)}_3R'{}^{(1)}\\
& =\wt{R}{}^{(2)}u\wt{R}{}^{(1)}_1R^{(2)}R'{}^{(2)}
S(R^{(1)})S(\wt{R}{}^{(1)}_2)\wt{R}{}^{(1)}_3R'{}^{(1)}\\
& =\wt{R}{}^{(2)}u\wt{R}{}^{(1)}R^{(2)}R'{}^{(2)}
S(R^{(1)})R'{}^{(1)}\\
& =\wt{R}{}^{(2)}u\wt{R}{}^{(1)}\\
&=\wt{R}{}^{(2)}S(\ov{R}^{(2)})\ov{R}^{(1)}\wt{R}{}^{(1)}\\
&=1.
\endalign
$$
Hence, $\wt{R}{}^{(2)}S^2(\wt{R}{}^{(1)})$ is a left inverse of
$u$. To prove that it is also  a right inverse, we note that by
condition~(3) of Definition~1.1 and Proposition~1.3,
$$
\allowdisplaybreaks \align S^2(\wt{R}{}^{(2)})uS^2(\wt{R}{}^{(1)})
&= uR^{(1)}\wt{R}{}^{(2)}_2
R'{}^{(1)}S(R'{}^{(2)})S(\wt{R}{}^{(2)}_1)R^{(2)}\wt{R}{}^{(2)}_3
S^2(\wt{R}{}^{(1)})\\
&=uR^{(1)}R'{}^{(1)}\wt{R}{}^{(2)}_1S(\wt{R}{}^{(2)}_2)
S(R'{}^{(2)})R^{(2)}\wt{R}{}^{(2)}_3S^2(\wt{R}{}^{(1)})\\
&=uR^{(1)}R'{}^{(1)}S(R'{}^{(2)})R^{(2)}\wt{R}{}^{(2)}S^2(\wt{R}{}^{(1)})\\
&=u\wt{R}{}^{(2)}S^2(\wt{R}{}^{(1)}).
\endalign
$$
So, $u\wt{R}{}^{(2)}S^2(\wt{R}{}^{(1)})= S^2(\wt{R}{}^{(2)})u
S^2(\wt{R}{}^{(1)})=\wt{R}{}^{(2)}u\wt{R}{}^{(1)}=1$, as we
want.\qed

\enddemo

\proclaim{Proposition 3.3} Let $(H,R)$ be a semiquasitriangular
Hopf algebra. The Drinfeld element $u$ satisfies
$$
\gather
\epsilon(u)=1, \quad \Delta(u)=(R_{21}R)^{-1}(u\ot u)=(u\ot u)(R_{21}R)^{-1},\\
\Delta(S(u))=(R_{21}R)^{-1}(S(u)\ot S(u))=(S(u)\ot
S(u))(R_{21}R)^{-1}
\intertext{and}
\Delta(uS(u))=(R_{21}R)^{-2}(uS(u)\ot uS(u))=(uS(u)\ot
uS(u))(R_{21}R)^{-2}.
\endgather
$$
\endproclaim

\demo{Proof} From Proposition~1.3 it is immediate that
$\epsilon(u)=1$. By Proposition~1.3 and item~(3) of
Definition~1.1, we have
$$
\align R_{32}R_{23}\bigl(R^{(1)}\ot & \Delta(S(R^{(2)})) \bigr) =
     R_{32}R_{23}\bigl(R^{(1)} \ot (S\ot S)\Delta^{\cop}(R^{(2)}) \bigr)\\
&=R_{32}(H\ot S\ot S)(R_{23})(H\ot S\ot S)\bigl( R^{(1)} \ot \Delta^{\cop}(R^{(2)}) \bigr)\\
&=R_{32}(H\ot S\ot S)\bigl(( R^{(1)}\ot  \Delta^{\op}(R^{(2)}))R_{23}\bigr)\\
&=R_{32}(H\ot S\ot S)\bigl(R_{23}( R^{(1)}\ot  \Delta(R^{(2)}))\bigr)\\
&=R_{32}(H\ot S\ot S)\bigl( R^{(1)}\ot  \Delta(R^{(2)})\bigr)R_{23}\\
&=(H\ot S\ot S)\bigl( (R^{(1)}\ot  \Delta(R^{(2)}))R_{32}\bigr)R_{23}\\
&=(H\ot S\ot S)\bigl( R_{32}(R^{(1)}\ot  \Delta^{\cop}(R^{(2)}))\bigr)R_{23}\\
&=\bigl(R^{(1)}\ot \Delta(S(R^{(2)})) \bigr)R_{32}R_{23}.
\endalign
$$
Hence, $R_{21}R\Delta(u) =R_{21}R \Delta(S(R^{(2)}))S(R^{(1)}) =
\Delta(S(R^{(2)}))R_{21}R S(R^{(1)})$. Similarly,
$\Delta(u)R_{21}R = \Delta(S(R^{(2)}))R_{21}R S(R^{(1)})$. On the
other hand arguing as in \cite{Ka, Proposition~VIII.4.5}, can be
check that this last expression equals $u\ot u$. This gives the
formula for $\Delta(u)$. Now it is easy to check the formulas for
$\Delta(S(u))$ and $\Delta(uS(u))$.\qed

\enddemo

\proclaim{Proposition 3.4} Let $(H,R)$ be a semiquasitriangular
Hopf algebra. The elements $u$, $S(u)$ and $uS(u)$ are
coinvariants for the coaction $\nu$.
\endproclaim

\demo{Proof} Since, by Proposition~2.3, $\nu$ is multiplicative,
it suffices to prove the assertion for $u$ and $S(u)$. We consider
the first case we leave the second one to the reader. By items~(3)
and (4) of Proposition~1.2 and Proposition~1.3,
$$
\align
\nu(u)&=R^{(2)}S(\ov{R}^{(2)}_2)\ov{R}^{(1)}_2R'{}^{(2)}\ot S(\ov{R}^{(1)}_1)
       S^2(\ov{R}^{(2)}_3)S(R^{(1)})S(\ov{R}^{(2)}_1)\ov{R}^{(1)}_3
       R'{}^{(1)}\\
&= R^{(2)}S(\ov{R}^{(2)}_2)R'{}^{(2)}\ov{R}^{(1)}_3\ot S(\ov{R}^{(1)}_1)
       S^2(\ov{R}^{(2)}_3)S(R^{(1)})S(\ov{R}^{(2)}_1)R'{}^{(1)}\ov{R}^{(1)}_2\\
&=
R^{(2)}S(R'{}^{(2)}\ov{R}^{(2)}_2)\ov{R}^{(1)}_3\ot S(\ov{R}^{(1)}_1)
       S^2(\ov{R}^{(2)}_3)S(R^{(1)})S(R'{}^{(1)}\ov{R}^{(2)}_1)\ov{R}^{(1)}_2\\
&=
R^{(2)}S(R'{}^{(2)})S(\ov{R}^{(2)}_1)\ov{R}^{(1)}_3\ot S(\ov{R}^{(1)}_1)
       S^2(\ov{R}^{(2)}_3)S(R^{(1)})S(R'{}^{(1)})S(\ov{R}^{(2)}_2)\ov{R}^{(1)}_2\\
&= S(\ov{R}^{(2)}_1)\ov{R}^{(1)}_3\ot S(\ov{R}^{(1)}_1)
       S^2(\ov{R}^{(2)}_3)S(\ov{R}^{(2)}_2)\ov{R}^{(1)}_2\\
&= S(\ov{R}^{(2)})\ov{R}^{(1)}_3\ot S(\ov{R}^{(1)}_1)\ov{R}^{(1)}_2\\
&=S(\ov{R}^{(2)})\ov{R}^{(1)}\ot 1\\
&=u\ot 1,
\endalign
$$
as we want.\qed
\enddemo

\proclaim{Corollary 3.5} For each semiquasitriangular
Hopf algebra $(H,R)$, it is true that $S^2(u)=u$, $S^2(u^{-1})=u^{-1}$ and
$uS(u)=S(u)u$.
\endproclaim

\demo{Proof} By Proposition~3.2 it suffices to prove that $T(u)=u$
and $T(S(u))=S(u)$. These assertions follow immediately from
Proposition~3.4, since $T = \mu\xcirc\nu$.\qed
\enddemo

\proclaim{Proposition 3.6} Let $(H,R)$ be a semiquasitriangular
Hopf algebra, and let $\wt{\nu}$ be the left coaction of $H$
defined by $\wt{\nu}:=(S\ot S)\xcirc \nu\xcirc S^{-1}$. The
following facts are equivalents:

\smallskip

\roster

\item $S(u)u \in \Z(H)$,

\smallskip

\item $S\xcirc T = T\xcirc S$,

\smallskip

\item $\mu\xcirc \nu=\mu \xcirc \wt{\nu}$.

\endroster

\endproclaim

\demo{Proof} Let $h\in H$. By Proposition~3.2, we have
$$
S(u)^{-1}S(T(h))S(u) S(S^2(h))=S^2(S(h)) = uT(S(h))u^{-1}.
$$
That is $S(T(h))S(u)u=S(u)uT(S(h))$. Since $S$ and $T$ are
bijective maps and $S(u)u$ is invertible, this implies that the
items~(1) and (2) are equivalents. It remains to  prove that (2)
$\Leftrightarrow$ (3). To do this it suffices to note that
$T=\mu\xcirc \nu$ and that
$$
S\xcirc T \xcirc S^{-1}= S\xcirc \mu \xcirc \nu \xcirc S^{-1} =
S\xcirc \mu^{\op} \xcirc \nu \xcirc S^{-1} = \mu \xcirc (S\ot S)
\xcirc \nu \xcirc S^{-1} = \mu\xcirc \wt{\nu},
$$
where the second equality follows from the fact that $\Ima \nu
\subseteq H\ot \Z(H)$.\qed

\enddemo

\enddocument